\documentclass[11pt]{article}
\usepackage{enumerate}
\usepackage{amssymb,a4wide,latexsym,makeidx,epsfig,fleqn}
\usepackage{amsthm}
\usepackage{amsmath}
\allowdisplaybreaks[4]
\usepackage{enumerate}
\usepackage{graphicx}
\usepackage{color}
\usepackage{epstopdf}
\newtheorem{theorem}{Theorem}[section]

\newtheorem{lemma}[theorem]{Lemma}

\begin{document}
\textwidth 150mm \textheight 230mm
\setlength{\topmargin}{-15mm}
\title{Scattering number and $\tau$-toughness in graphs involving $A_\alpha$-spectral radius
\footnote{This work is supported by the National Natural Science Foundations of China (No. 12371348, 12201258), the Postgraduate Research \& Practice Innovation Program of Jiangsu Normal University (No. 2025XKT0632, 2025XKT0633) .}}
\author{{ Caili Jia, Yong Lu\footnote{Corresponding author.}}\\
{\small  School of Mathematics and Statistics, Jiangsu Normal University,}\\ {\small  Xuzhou, Jiangsu 221116,
People's Republic
of China.}\\
{\small E-mails: jiacaili0309@163.com, luyong@jsnu.edu.cn}}

\date{}
\maketitle
\begin{center}
\begin{minipage}{120mm}
\vskip 0.3cm
\begin{center}
{\small {\bf Abstract}}
\end{center}
{\small
The scattering number $s(G)$ of graph $G=(V,E)$ is defined as $s(G)$=max\big\{$c(G-S)-|S|$\big\}, where the maximum is taken over all proper subsets $S\subseteq V(G)$, and $c(G-S)$ denotes the number of components of $G-S$.  In 1988, Enomoto introduced a variation of toughness $\tau(G)$ of a graph $G$, which is defined by $\tau(G)$=min\big\{$\frac{|S|}{c(G-S)-1}$, $S\subseteq V(G)$ and $c(G-S)>1$\big\}. Both the scattering number and toughness are used to characterize the invulnerability or stability of a graph, i.e., the ability of a graph to remain connected after vertices or edges are removed. The smaller the value of $s(G)$ (or the larger the value of $\tau(G)$), the stronger the connectivity of a graph $G$. The $A_{\alpha}$-spectral radius of $G$ is denoted by $\rho_{\alpha}(G)$. Using typical $A_{\alpha}$-spectral techniques and structural analysis, we present a sufficient condition such that $s(G)\leq 1$. This result generalizes the result of Chen, Li and Xu [Graphs Comb. 41 (2025)]. Furthermore, we establish a sufficient condition with respect to the $A_{\alpha}$-spectral radius for a graph to be $\tau$-tough. When  $\alpha=\frac{1}{2}$, our result reduces to that of  Chen, Li and Xu [Graphs Comb. 41 (2025)].
 }
\end{minipage}
\end{center}

\section{Introduction }
\hspace{1.3em}
Throughout this paper, we consider only finite, undirected and simpe graphs. Let $G=(V(G),E(G))$ be a graph, where $V(G)$ is the vertex set and $E(G)$ is the edge set.
The \emph{order} and \emph{size} of $G$ are denoted by $|V(G)|=n$ and $|E(G)|=m$, respectively. Denote by $\delta(G)$ ($\delta$ for short) and $\Delta(G)$ ($\Delta$ for short) the minimum degree and the maximum degree of $G$, respectively.
For a vertex subset $S$ of $G$, we denote by $G-S$ and $G[S]$ the subgraph of $G$ obtained from $G$ by deleting the vertices in $S$ together with their incident edges and the subgraph of $G$ induced by $S$, respectively.
A subgraph of a graph $G$ is \emph{spanning} if the subgraph covers all vertices of $G$.
The number of components of $G$ is denoted by $c(G)$.

Let $G$ be a graph of order $n$, and let the \emph{adjacency matrix} of $G$ be defined as $A(G)=(a_{ij})_{n\times n}$, where $a_{ij}=1$ if $v_{i}v_{j}\in E(G)$, and $a_{ij}=0$ otherwise.
The \emph{degree diagonal matrix} is the diagonal matrix of vertex degrees of $G$, denoted by $D(G)$. The \emph{Laplacian matrix} $L(G)$ and \emph{signless Laplacian matrix} $Q(G)$ of $G$ is defined by $L(G)=D(G)-A(G)$ and $Q(G)=D(G)+A(G)$, respectively.
The eigenvalues of $A(G)$, $L(G)$ and $Q(G)$ are called the \emph{eigenvalues}, the \emph{Laplacian eigenvalues} and the \emph{signless Laplacian eigenvalues} of $G$, and denoted by $\lambda_{1}(G)\geq \lambda_{2}(G)\geq\cdots\geq\lambda_{n}(G)$, $\mu_{1}(G)\geq \mu_{2}(G)\geq\cdots\geq\mu_{n}(G)=0$ and $q_{1}(G)\geq q_{2}(G) \geq\cdots\geq q_{n}(G)$, respectively. The largest eigenvalues of $A(G)$, $L(G)$ and $Q(G)$ are also called the \emph{spectral radius}, \emph{Laplacian spectral radius} and the \emph{signless Laplacian spectral radius} of $G$, and denoted by $\lambda(G)$, $\mu(G)$ and $q(G)$, respectively. We sometimes write $\lambda_i$, $\mu_i$ and $q_i$ to instead of $\lambda_{i}(G), \mu_{i}(G)$ and $q_i(G)$ for $1\leq i \leq n$.
For any $\alpha\in[0,1]$, Nikiforov \cite{N} introduced the \emph{$A_{\alpha}$-matrix} of $G$ as $A_{\alpha}(G)=\alpha D(G)+(1-\alpha)A(G)$. One can see that $A_{\alpha}(G)=A(G)$ if $\alpha=0$, and $A_{\alpha}(G)=\frac{1}{2}Q(G)$ if $\alpha=\frac{1}{2}$. Hence, $A_{\alpha}(G)$ generalizes both the adjacency matrix and the signless Laplacian matrix of $G$. The eigenvalues of $A_{\alpha}(G)$ are called the \emph{$A_{\alpha}$-eigenvalues} of $G$, and the largest of them, denoted by $\rho_{\alpha}(G)$, is called the \emph{$A_{\alpha}$-spectral radius} of $G$. More interesting spectral properties of $A_{\alpha}(G)$ can be found in \cite{BFO,LL,LHX,LLX,N,NR}.

Let $t$ be a positive real number and a connected graph $G$ is \emph{$t$-tough} if $tc(G-S)\leq|S|$ for every vertex cut $S$ of $V(G)$.
The \emph{toughness} of graph $G$, denoted by $t(G)$, is the largest value of $t$ for which $G$ is $t$-tough (taking $t(K_{n})=\infty$ where $K_{n}$ is a complete graph of order $n$).
Thus, $t(G)=\min\{\frac{|S|}{c(G-S)}:S\subseteq V(G),c(G-S)>1\}$.
The concept of toughness initially proposed by Chv\'{a}tal \cite{C} in 1973, which serves as a simple way to measure how tightly various pieces of a graph hold together. The toughness is related to many other important properties of a graph, such as the existence of factors \cite{B,CL,E,EH,EJKS,LQSX}, cycles \cite{FN} and spanning trees \cite{BBHV}. For more extensive work on toughness, one can see \cite{BBS,CLW,LFS,MHY,S,SS}.
In order to better investigate the existence of factors in a graph, Enomoto \cite{EH} introduced a slight variation of toughness in
1998. A non-complete graph $G$ is \emph{$\tau$-tough} if $|S|\geq\tau(c(G-S)-1)$ for every proper subset $S\subseteq V(G)$ with $c(G-S)>1$. The
\emph{variation of toughness} $\tau(G)$ of $G$ is the maximum $\tau$ for which $G$ is $\tau$-tough. Thus, $\tau(G)=\min\{\frac{|S|}{c(G-S)-1}:S\subseteq V(G),c(G-S)>1\}$ in which the minimum is taken over all proper sets $S\subseteq V(G)$.
For two vertex-disjoint graphs $G_{1}$ and $G_{2}$, we use $G_{1}\cup G_{2}$ to denote the \emph{disjoint union} of $G_{1}$ and $G_{2}$.
The \emph{join} $G_{1}\vee G_{2}$ is the graph obtained from $G_{1}\cup G_{2}$ by adding all possible edges between $V(G_{1})$ and $V(G_{2})$.

In \cite{C}, Chv\'{a}tal conjectured that every $k$-tough graph on $n\geq k+1$ vertices and $kn$ even contains a $k$-factor. Enomoto et al. \cite{EJKS} gave a decisive answer to Chv\'{a}tal's conjecture. In \cite{E}, Enomoto strengthened the result obtained by  Enomoto et al. \cite{EJKS}. In 1998, this result was first improved by Enomoto \cite{EH} by using the definition of $\tau(G)$. Enomoto and Hagita \cite{EHH} were able to generalize the results in \cite{EH} and strengthened the results in \cite{EJKS} with a sufficiently large number of vertices.

The \emph{scattering number} of a graph $G$ was defined by Jung \cite{J} as $s(G)$=max\big\{$c(G-S)-|S|: S\subseteq V(G),c(G-S)>1$\big\}.
It is clear that the concept of a scattering number is similar to that of toughness in some sense.
Hendry \cite{G} proposed $s(G)\leq0$ is one necessary condition for Hamiltonian since $s(G)\leq0$ is equivalent to $t(G)\geq1$.
For any graph $G$ the condition $s(G)\leq 1$ is equivalent to the condition $\tau(G)\geq1$.
Both the scattering number and toughness are used to characterize the invulnerability or stability of a graph, i.e., the ability of the graph to remain connected after vertices or edges are removed.
The parameter of scattering number has a strong background of applications in measuring network vulnerability \cite{KK}. Unlike the connectivity measures, it takes account not only the difficulty to destroy a network, but also what remains after the network is destroyed. Therefore, this parameter attracts much attention from researchers, one can see \cite{C,F,G,JRZ,KK,ZW}.
Jamrozik, Kalinowski and Skupien \cite{JRZ} studied the small maximal non-Hamiltonian graphs by using a scattering number.
In \cite{G}, Hendry used the concept of a scattering number to study extremal non-Hamiltonian graphs. It is found that the concept of a scattering number is more convenient than the notion of toughness for describing maximal and extremal non-Hamiltonian graphs. Zhang, Li and Han \cite{ZLH} demonstrated that calculating a graph's scattering number is an NP-complete problem. Consequently, exploring the bounds for the scattering number of graphs becomes a matter of significant interest.

Recently, Li, Shi and Gu \cite{LSG} investigated the relationship between the scattering number of a regular graph and its eigenvalues, and derived an upper bound for $s(G)$ in the case of regular graphs $G$.  Moreover, Gu and Liu \cite{GL} proved that for any graph $G$ of order $n$, $s(G)\leq \max\Big\{0, \frac{(\mu_{1}-2\mu_{n-1})n+2\mu_{n-1}}{\mu_{1}}\Big\}$. On the other hand, by applying classical spectral methods, Fan, Lin and Lu \cite{FLL} derived a tight lower bound for $\lambda(G)$ that ensures a graph G satisfying $t(G)\geq 1$ (or $s(G)\leq 0$). Inspired by their findings, Chen, Li and Shiu \cite{CLS} went on to establish a tight lower bound for $q(G)$ to guarantee that a graph $G$ meets $s(G)\leq 0$. Noted that $s(G)\leq1$ constitutes a well-known necessary condition for a graph to be traceable. It is natural to pose the question: ``Do there exist any $A_{\alpha}$-spectral conditions that can ensure a graph $G$ has $s(G)\leq1$?" Chen, Li and Xu \cite{CLX} established a sufficient condition involving the spectral radius for a graph $G$ with minimum degree $\delta$ such that $s(G)\leq1$, as well as a condition involving Laplacian eigenvalues for a graph with $s(G)\leq1$.

Focus on this problem, we further develop the following tight lower bound for $\rho_{\alpha}(G)$ to guarantee that a graph $G$ satisfies $s(G)\leq1$.

\noindent\begin{theorem}\label{th:1.1.}
Let $G$ be a connected graph of order $n\geq$ max$\{4\delta+2,\delta^{3}+\delta\}$ with minimum degree $\delta$.
If $$\rho_{\alpha}(G)\geq\rho_{\alpha}(K_{\delta}\vee(K_{n-2\delta-1}\cup{(\delta+1)} K_{1})),$$ then $s(G)\leq 1$ unless $G\cong K_{\delta}\vee(K_{n-2\delta-1}\cup{(\delta+1)} K_{1})$.
\end{theorem}

In particular, if $\alpha=0$, then we can generalize the theorem as follows which confirmed by  Chen, Li and Xu in \cite{CLX}.

\textbf{[Theorem 1.1 in \cite{CLX}]}
Let $G$ be a connected graph of order $n\geq$max$\{\delta^3+\delta, 8\delta+4$\} with minimum degree $\delta$. If
\begin{center}
$\lambda(G)\geq \lambda(K_{\delta}\vee(K_{n-2\delta-1}\cup{(\delta+1)} K_{1}))$,
\end{center}
then $s(G)\leq 1$ unless $G\cong K_{\delta}\vee(K_{n-2\delta-1}\cup{(\delta+1)} K_{1})$.

It was showed in \cite{GWC} that the feasibility of network data  transmission is equivalent to the existence of fractional flow in the network, and after modeling with graphs, it can be characterized by the existence of factors or fractional factors in various settings. The previous works can fully explain that during the network design and construction stage, if the network graph satisfies a certain toughness condition, theoretically, it can simultaneously guarantee the network's robustness and the feasibility of data transmission in the network \cite{GWC}.
Recently, Gu and Liu \cite{GL} showed that for any connected graph $G$, $\tau(G)>\frac{\mu_{n-1}}{\mu_{1}{\mu_{n-1}}}$. It is natural to ask whether there are some sufficient conditions for a graph $G$ to be $\tau$-tough, or there exists a tight lower bound for $\tau(G)$? By incorporating the variation of toughness and the spectral radius, Chen, Fan and Lin \cite{CFL} provided conditions involving the spectral radius for a graph to be $\tau$-tough ($\tau\geq 2$ is an integer) and to be $\tau$-tough ($\frac{1}{\tau}$ is a positive integer), respectively. In \cite{CLX},  Chen, Li and Xu presented two lower bounds on the size to guarantee a graph $G$ to be $\tau$-tough and
further provided the signless Laplacian spectral conditions for a graph to be $\tau$-tough. Recently, Zhou et al. \cite{ZZZL} presented an $A_{\alpha}$-specatral radius condition for a graph to be $t$-tough. In the paper, we prove two lower bounds on the $A_{\alpha}$-spectral radius to guarantee a graph $G$ to be $\tau$-graph and construct the corresponding extremal graphs to show all these bounds are best possible.

\noindent\begin{theorem}\label{th:1.2.}
Suppose that $G$ is a connected graph of order $n$. Then the following statements hold.
\begin{enumerate}[(1)]
\item
Let $\tau\geq 2$ be an integer, $\alpha\in[\frac{1}{2},\frac{3}{4})$ and $n\geq\max\{{4\tau^2+5\tau+1,\frac{8\tau(1-\alpha)-2\alpha+1}{3-4\alpha}}\}$. If $\rho_{\alpha}(G)\geq\rho_{\alpha}(K_{\tau-1}\vee(K_{n-\tau}\cup K_{1}))$, then  $G$  is a $\tau$-tough graph unless $G\cong K_{\tau-1}\vee(K_{n-\tau}\cup K_{1})$.

\item
Let $\frac{1}{\tau}\geq 2$ be a positive integer,
$\alpha\in[\frac{1}{2},\frac{3+\tau}{4+2\tau})$
and
$n\geq\max\{ 2\tau^2+5\tau+\frac{2}{\tau}+8, T$\}, where $T=\frac{(5-6\alpha)\tau^2+(13-14\alpha)\tau+4(1-\alpha)}{(1-2\alpha)\tau^2+(3-4\alpha)\tau}$.
If $\rho_{\alpha}(G)\geq\rho_{\alpha}(K_{1}\vee(K_{n-\frac{1}{\tau}-2}\cup(\frac{1}{\tau}+1) K_{1}))$, then  $G$  is a $\tau$-tough graph unless $G\cong K_{1}\vee(K_{n-\frac{1}{\tau}-2}\cup(\frac{1}{\tau}+1) K_{1})$.
\end{enumerate}
\end{theorem}

Since $A_\frac{1}{2}(G)=\frac{1}{2}Q(G)$, it follows that $\rho_\frac{1}{2}(G)=\frac{1}{2}q(G)$. In particular, if $\alpha=\frac{1}{2}$, then we can derive the theorem as follows which confirmed by  Chen, Li and Xu in \cite{CLX}.

\textbf{[Theorem 1.7 in \cite{CLX}]}
Let $G$ be a connected graph of order $n$. The following statements holds.
\begin{enumerate}[(1)]
\item
Let $\tau\geq 2$ be an integer and $n\geq(4\tau+1)(\tau+1)$. If $q(G)\geq q(K_{\tau-1}\vee(K_{n-\tau}\cup K_{1}))$, then  $G$  is a $\tau$-tough graph unless $G\cong K_{\tau-1}\vee(K_{n-\tau}\cup K_{1})$.

\item
Let $\frac{1}{\tau}\geq 1$ be a positive integer and $n\geq(2\tau+5)\tau+\frac{2}{\tau}+8$. If $q(G)\geq q(K_{1}\vee(K_{n-\frac{1}{\tau}-2}\cup(\frac{1}{\tau}+1) K_{1}))$, then  $G$  is a $\tau$-tough graph unless $G\cong K_{1}\vee(K_{n-\frac{1}{\tau}-2}\cup(\frac{1}{\tau}+1) K_{1})$.
\end{enumerate}


\section{Proof of Theorem \ref{th:1.1.} }
\hspace{1.3em}
In this section, we introduce some useful lemmas, which play an important role in the proof of Theorem \ref{th:1.1.}.\\

\noindent\begin{lemma}\label{le:2.1.}\cite{N}
Let $K_{n}$ be a complete graph of order $n$. Then
\begin{align*}
\rho_{\alpha}(K_{n})=n-1.
\end{align*}
\end{lemma}

\noindent\begin{lemma}\label{le:2.2.}\cite{N}
If $G$ is a connected graph, and $H$ is a proper subgraph of $G$. Then
\begin{align*}
\rho_{\alpha}(G)>\rho_{\alpha}(H).
\end{align*}
\end{lemma}

Let $M$ be a real $n\times n$ matrix, and let $X=\{1,2,\ldots,n\}$. Given a partition $\Pi=\{X_1, X_2,\ldots, X_k\}$ with $X=X_1\cup X_2\cup\cdots\cup X_k$, the matrix $M$ can be correspondingly partitioned as
\begin{align*}
M=\left(
\begin{array}{ccccccccc}
M_{11} & M_{12} & \cdots & M_{1k} \\
M_{21} & M_{22} & \cdots & M_{2k} \\
\vdots & \vdots & \ddots &\vdots \\
M_{k1} & M_{k2} & \cdots & M_{kk} \\
\end{array}
\right),
\end{align*}
The \emph{quotient matrix} of $M$ with respect to $\Pi$ is defined as the $k\times k$ matrix $B_{\Pi}=(b_{i,j})_{i,j=1}^{k}$
where $b_{i,j}$ is the average value of all row sums of $M_{i,j}$. The partition $\Pi$ is called \emph{equitable} if each block $M_{i,j}$ of $M$ has constant row sum $b_{i,j}$. Also, we say that the quotient matrix $B_{\Pi}$ is \emph{equitable} if $\Pi$ is an equitable partition of $M$.

\noindent\begin{lemma}\label{le:2.3.}\cite{YYSX}
Let $M$ be a nonnegative matrix, and let $B$ be an equitable quotient matrix of $M$. Then the eigenvalues of $B$ are also eigenvalues
of $M$, and
\begin{align*}
\lambda(M)=\lambda(B).
\end{align*}
\end{lemma}

\noindent\begin{lemma}\label{le:2.4.}\cite{GR}
If $M_1$ and $M_2$ are two nonnegative $n\times n$ matrices such that $M_1-M_2$ is nonnegative, then
\begin{align*}
\lambda(M_1)\geq\lambda(M_2),
\end{align*}
where $\lambda(M_{i})$ is the spectral radius of $M_{i}$ for $i =1, 2$.
\end{lemma}

\noindent\begin{lemma}\label{le:2.5.}\cite{ZHW}
Let $\alpha\in[0,1)$, and let $n_{1}\geq n_{2}\geq\cdots\geq n_{t}$ be positive integers with $n=\sum\limits_{i=1}^{t}n_{i}+s$ and $n_{1}\leq n-s-t+1$. Then
\begin{align*}
\rho_{\alpha}(K_{s}\vee(K_{n_{1}}\cup K_{n_{2}}\cup\cdots\cup K_{n_{t}}))\leq\rho_{\alpha}(K_{s}\vee(K_{n-s-t+1}\cup(t-1) K_{1})),
\end{align*}
where the equality holds if and only if $(n_1,n_2,\ldots,n_t)=(n-s-t+1,1,\ldots,1)$.
\end{lemma}

For the proof of Theorem \ref{th:1.1.}, we need the following critical lemma, which generalizes Lemma \ref{le:2.5.}.\\

\noindent\begin{lemma}\label{le:2.6.}
Let $\alpha\in[0,1)$, and let $n_{1}\geq n_{2}\geq\cdots\geq n_{t}\geq p$ be positive integers with $n=\sum\limits_{i=1}^{t}n_{i}+s$ and $n_{1}\leq n-s-p(t-1)$. Then
\begin{align*}
\rho_{\alpha}(K_{s}\vee(K_{n_{1}}\cup K_{n_{2}}\cup\cdots\cup K_{n_{t}}))\leq\rho_{\alpha}(K_{s}\vee(K_{n-s-p(t-1)}\cup(t-1) K_{p})),
\end{align*}
where the equality holds if and only if $(n_1,n_2,\ldots,n_t)=(n-s-p(t-1),p,\ldots,p)$.\\
\end{lemma}

\noindent\textbf{Proof.}
It suffices to prove
\begin{align*}
\rho_{\alpha}(K_{s}\vee(K_{n_1}\cup K_{n_2}\cup\cdots\cup K_{n_t}))<\rho_{\alpha}(K_{s}\vee(K_{n_{1}+1}\cup K_{n_2}\cup\cdots\cup K_{n_{j}-1}\cup\cdots\cup K_{n_t}))
\end{align*}
where $n_{j}\geq p$ for $j\in\{2,\ldots,t\}$. Without loss of generality, we take $j=t$, and for other cases the proof is similar.

Let $G_{1}=K_{s}\vee(K_{n_1}\cup K_{n_2}\cup\cdots\cup K_{n_t})$ and $G_{2}=K_{s}\vee(K_{n_{1}+1}\cup K_{n_2}\cup\cdots\cup K_{n_{t}-1})$.
Observe that the equitable quotient matrix of $A_{\alpha}(G_{1})$ equals
\begin{align*}
B_{1}=
\bordermatrix{%
& s & n_{1}& n_{2} & \cdots & n_{t}\cr
s &n\alpha-s\alpha+s-1 & n_{1}(1-\alpha)& n_{2}(1-\alpha) & \cdots & n_{t}(1-\alpha) \cr
n_{1} & s(1-\alpha) & s\alpha+n_{1}-1 & 0 & \cdots & 0\cr
n_{2} & s(1-\alpha) & 0 & s\alpha+n_{2}-1 & \cdots & 0\cr
\vdots & \vdots & \vdots & \vdots & \ddots &\vdots \cr
n_{t} & s(1-\alpha) & 0 & 0 & \cdots & s\alpha+n_{t} -1\cr
}.
\end{align*}
Then the characteristic polynomial of $B_1$ is equal to
\begin{align*}
f_{B_{1}}(x)=&(x-n\alpha+s\alpha-s+1)\cdot\prod_{j=1}^{t}(x-s\alpha-n_{j}+1)
\\&+s(1-\alpha)^2\cdot\sum_{i=1}^{t}(-1)^{i}n_{i}\prod_{j\neq i}(x-s\alpha-n_{j}+1)
\\=&\varphi(s,n_1,n_2,\ldots,n_t;x).
\end{align*}
By Lemma\ref{le:2.3.}, $\rho_\alpha(G_1)$ equals the largest root of $f_{B_{1}}(x)= 0$. Then, by Lemma\ref{le:2.4.}, we obtain
\begin{align}
\rho_{\alpha}(G_{1})>n\alpha-s\alpha+s-1.
\end{align}
Since $G_1$ contains $K_{n_{1}+s}$ as a proper subgraph, we have
\begin{align}
\rho_{\alpha}(G_{1})>\rho_{\alpha}(K_{n_{1}+s})=n_{1}+s-1.
\end{align}
Also, note that $A_{\alpha}(G_{2})$ has the equitable quotient matrix
\begin{align*}
B_{2}=
\bordermatrix{%
& s & n_{1}+1& n_{2} & \cdots & n_{t}-1\cr
s &n\alpha-s\alpha+s-1 & (n_{1}+1)(1-\alpha)& n_{2}(1-\alpha) & \cdots & (n_{t}-1)(1-\alpha) \cr
n_{1}+1 & s(1-\alpha) & s\alpha+n_{1} & 0 & \cdots & 0\cr
n_{2} & s(1-\alpha) & 0 & s\alpha+n_{2}-1 & \cdots & 0\cr
\vdots & \vdots & \vdots & \vdots & \ddots &\vdots \cr
n_{t}-1 & s(1-\alpha) & 0 & 0 & \cdots & s\alpha+n_{t} -2\cr
}.
\end{align*}
The characteristic polynomial of $B_2$ is equal to
\begin{align*}
f_{B_{2}}(x)=\varphi(s,n_{1}+1,n_2,\ldots,n_{t}-1;x).
\end{align*}
As above, we see that $\rho_{\alpha}(G_2)$ is equal to the largest root of $f_{B_2}(x)=0$. By a computation, we obtain
\begin{align*}
f_{B_{2}}(x)-f_{B_{1}}(x)&=\varphi(s,n_{1}+1,n_2,\ldots,n_{t}-1;x)-\varphi(s,n_1,n_2,\ldots,n_t;x)
\\&=b_{0}(x)+s(1-\alpha)^2\sum_{i=1}^{t}b_{i}(x),
\end{align*}
where
\begin{align*}
&b_{0}(x)=-(n_{1}-n_{t}+1)(x-n\alpha+s\alpha-s+1)\prod_{j=2}^{t-1}(x-s\alpha-n_{j}+1),\\
&b_{1}(x)=-(x-s\alpha+n_{1}-n_{t}+2)\prod_{j=2}^{t-1}(x-s\alpha-n_{j}+1),\\
&b_{i}(x)=(-1)^{i+1}n_{i}(n_{1}-n_{t}+1)\prod_{j=2,j\neq i}^{t-1}(x-s\alpha-n_{j}+1) \enspace for \enspace i=2,\ldots,t-1,\\
&b_{t}(x)=(-1)^{t+1}(x-s\alpha+n_{1}-n_{t})\prod_{j=2}^{t-1}(x-s\alpha-n_{j}+1).
\end{align*}
Combining $(1)$, $(2)$ and $n_1\geq n_2\geq\cdots\geq n_t$, we get $b_{0}(\rho_{\alpha}(G_1))\leq 0$.\\

\textbf{Case 1.} $t$ is odd.

Combining $(2)$, we obtain
\begin{align*}
&b_{1}(\rho_{\alpha}(G_1))+b_{t}(\rho_{\alpha}(G_1))=-2(n_{1}-n_{t}+1)\prod_{j=2}^{t-1}(\rho_{\alpha}(G_1)-s\alpha-n_{j}+1)<0,\\
&b_{t-1}(\rho_{\alpha}(G_1))=-n_{t-1}(n_{1}-n_{t}+1)\prod_{j=2,j\neq i}^{t-1}(\rho_{\alpha}(G_1)-s\alpha-n_{j}+1)<0,\\
&b_{i}(\rho_{\alpha}(G_1))+b_{i+1}(\rho_{\alpha}(G_1))=-(n_{i}-n_{i+1})(n_{1}-n_{t}+1)(\rho_{\alpha}(G_1)-s\alpha+1)\cdot
\\ & \prod_{j=2,j\neq i,i+1}^{t-1}(\rho_{\alpha}(G_1)-s\alpha-n_{j}+1)\leq 0,
\end{align*}
for all even $i$ with $2\leq i\leq t-3$.\\

\textbf{Case 2.} $t$ is even.

Combining $(2)$, we obtain
\begin{align*}
&b_{1}(\rho_{\alpha}(G_1))+b_{t}(\rho_{\alpha}(G_1))=-2(n_{1}-n_{t}+1)\prod_{j=2}^{t-1}(\rho_{\alpha}(G_1)-s\alpha-n_{j}+1)<0,\\
&b_{i}(\rho_{\alpha}(G_1))+b_{i+1}(\rho_{\alpha}(G_1))=-(n_{i}-n_{i+1})(n_{1}-n_{t}+1)(\rho_{\alpha}(G_1)-s\alpha+1)\cdot
\\ & \prod_{j=2,j\neq i,i+1}^{t-1}(\rho_{\alpha}(G_1)-s\alpha-n_{j}+1)\leq 0,
\end{align*}
for all even $i$ with $2\leq i\leq t-2$.
Since $f_{B_1}(\rho_{\alpha}(G_1))=0$, we deduce that
\begin{align*}
f_{B_2}(\rho_{\alpha}(G_1))=f_{B_2}(\rho_{\alpha}(G_1))-f_{B_1}(\rho_{\alpha}(G_1))< 0,
\end{align*}
which implies that $\rho_{\alpha}(G_1)<\rho_{\alpha}(G_2)$.

The proof is completed.
$\hfill\square$\\

Now we shall give a proof of Theorem \ref{th:1.1.}.\\

\noindent\textbf{Proof of Theorem \ref{th:1.1.}.}

Let $G$ be a graph satisfying the conditions in Theorem 1.1, and let $H=K_{\delta}\vee(K_{n-2\delta-1}\cup (\delta+1)K_1)$. Note that $s(H)=2$.
Suppose to the contrary, we assume that $s(G)>1$ with $G\neq H$, show that $\rho_{\alpha}(G)<\rho_\alpha(H)$.
Thus, there exists a nonempty vertex set $S\subseteq V(G)$ such that $c(G-S)-|S|>1$, i.e., $c(G-S)\geq|S|+2$.
Let $|S|=s$. We know $G$ is a spanning subgraph of $G_s^1=K_{s}\vee(K_{n_{1}}\cup K_{n_{2}}\cup\cdots\cup K_{n_{s+2}})$ for some positive integers $n_{1}\geq n_{2}\geq\cdots\geq n_{s+2}\geq1$ and $\sum\limits_{i=1}^{s+2}n_{i}=n-s$.
By Lemma \ref{le:2.2.}, we obtain
\begin{align}
\rho_{\alpha}(G)\leq\rho_\alpha(G_s^1)
\end{align}
with equality  if and only if $G\cong G_s^1$.
We now consider the following three cases.\\

\textbf{Case 1.} $s=\delta$.

Note that $G_s^1$ is a spanning subgraph of $H$. Thus
\begin{align*}
\rho_{\alpha}(G_s^1)\leq\rho_\alpha(H)
\end{align*}
with equality if and only if $G_s^1=H$. Combining this with ($3$), we conclude that
\begin{align*}
\rho_{\alpha}(G)\leq\rho_\alpha(H)
\end{align*}
with equality if and only if $G=H$. Since $G\neq H$, we have $\rho_{\alpha}(G)<\rho_\alpha(H)$, a contradiction.\\

\textbf{Case 2.} $s<\delta$.

Let $G_s^2=K_{s}\vee(K_{n-s-(\delta+1-s)(s+1)}\cup (s+1)K_{\delta+1-s})$. Recall that $G$ is a spanning subgraph of $G_s^1=K_{s}\vee(K_{n_{1}}\cup K_{n_{2}}\cup\cdots\cup K_{n_{s+2}})$, where $n_{1}\geq n_{2}\geq\cdots\geq n_{s+2}$ and $\sum\limits_{i=1}^{s+2}n_{i}=n-s$. Clearly, $n_{s+2}\geq \delta+1-s$ because the minimum degree of $G_s^1$ is at least $\delta$. By Lemma \ref{le:2.6.}, we have
\begin{align}
\rho_{\alpha}(G_s^1)\leq \rho_{\alpha}(G_s^2)
\end{align}
with equality if and only if $(n_1, n_2, \ldots, n_{s+2})=(n-s-(\delta+1-s), \delta+1-s, \ldots, \delta+1-s)$. Note that the $A_\alpha$-matrix of $G_s^2$ is as follows
\begin{align*}
\scalebox{0.8}{$
\bordermatrix{%
& n-s-(\delta+1-s)(s+1) & (s+1)(\delta+1-s) & s \cr
n-s-(\delta+1-s)(s+1)    & n-(\delta+2-s)(s+1)+\alpha s & 0 & (1-\alpha) s \cr
(s+1)(\delta+1-s) & 0 & (\delta-s)(s(1-\alpha)+1)+\alpha s & (1-\alpha)s \cr
s  & (1-\alpha)(n-s-(\delta+1-s)(s+1))& (1-\alpha)(\delta+1-s)(s+1) & \alpha(n-s)+s-1 \cr
}.$}
\end{align*}

Suppose $\rho_{\alpha}(G_s^2)\geq n-(\delta+1-s)(s+1)$. Let $\rho_{\alpha}(G_s^2)=\rho$ and let $\textbf{x}$ be the Perron vector of $A_{\alpha}(G_s^2)$ with respect of $\rho$. By symmetry, we take $x_u=x_1$ for all $u\in V(K_{n-s-(\delta+1-s)(s+1)})$, $x_v=x_2$ for all $v\in V((s+1)K_{\delta+1-s})$, and $x_w=x_3$ for all $w\in V(K_s)$. According to $A_\alpha(G_s^2)\textbf{x}=\rho\textbf{x}$, we have

\begin{small}
\begin{align}
&\rho x_{1}=(n-(\delta+2-s)(s+1)+\alpha s)x_{1}+(1-\alpha)sx_{3},\\
&\rho x_{2}=\big[(\delta-s)\big(s(1-\alpha)+1\big)+\alpha s\big]x_{2}+(1-\alpha)sx_{3},\\
&\rho x_{3}=(1-\alpha)\Big[\big(n-s-(\delta+1-s)(s+1)\big)x_{1}+(\delta+1-s)(s+1)x_{2}\Big]+(\alpha(n-s)+s-1)x_{3}.
\end{align}
\end{small}

From ($5$) and ($6$), we obtain

\begin{align}
&x_{1}=\frac{(1-\alpha)s}{\rho-(n-(\delta+2-s)(s+1)+\alpha s)}x_{3},\\
&x_{2}=\frac{(1-\alpha)s}{\rho-((\delta-s)(s(1-\alpha)+1)+\alpha s)}x_{3}.
\end{align}

Since $n\geq \delta^3+\delta$ and $\delta>s$,
\begin{align*}
\rho &\geq n-(\delta+1-s)(s+1)\nonumber
\\&\geq \delta^3+\delta-(\delta+1)(s+1)+s(s+1)\nonumber
\\&=\delta^3-(\delta+1)(s+1)+s(s+1)-1+\delta+1\nonumber
\\&>\delta^3-(\delta+1)^2+1+\delta+1\nonumber
\\&=(\delta^2-\delta-2)\delta+\delta+1\nonumber
\\&>\delta+1.\nonumber
\end{align*}

Putting ($8$) and ($9$) into ($7$), we get
\begin{align}
\rho+1=&s+\alpha(n-s)+\frac{(1-\alpha)^2s(n-s-(\delta+1-s)(s+1))}{\rho-(n-(\delta+1-s)(s+1)-1)+(1-\alpha)s}\nonumber
\\&+\frac{(1-\alpha)^2s(\delta+1-s)(s+1)}{\rho-((\delta-s)(s(1-\alpha)+1)+\alpha s)}\nonumber
\\ \leq& s+\alpha(n-s)+\frac{(1-\alpha)^2s(n-s-(\delta+1-s)(s+1))}{n-(\delta+1-s)(s+1)-(n-(\delta+1-s)(s+1)-1)+(1-\alpha)s}\nonumber
\\&+\frac{(1-\alpha)^2s(\delta+1-s)(s+1)}{n-(\delta+1-s)(s+1)-((\delta-s)(s(1-\alpha)+1)+\alpha s)}\nonumber
\\<&s+\alpha(n-s)+\frac{(1-\alpha)^2s(n-s-(\delta+1-s)(s+1))}{(1-\alpha)s+1}+\frac{(1-\alpha)^2s(\delta+1-s)(s+1)}{(1-\alpha)s+1}\nonumber
\\=&s+\alpha(n-s)+(1-\alpha)(n-s)-\frac{(1-\alpha)(n-s)}{(1-\alpha)s+1}\nonumber
\\=&n-\frac{n-s}{s+\frac{1}{1-\alpha}}\nonumber
\\ \leq& n-\frac{n-s}{s+1}
\\=&n-(\delta+1-s)(s+1)-\frac{n-s-(\delta+1-s)(s+1)^2}{s+1}
\\ \leq & n-(\delta+1-s)(s+1)-\frac{\delta^3+\delta-s-(\delta+1-s)(s+1)^2}{s+1}\nonumber.
\end{align}

Suppose $s=1$, $\rho\geq n-2\delta$. From ($10$), we have
\begin{center}
$\rho+1<n-\frac{n-1}{2}<n-\frac{\delta^3+\delta-1}{2}<n-2\delta\leq\rho$,
\end{center}
a contradiction.

Suppose $s=2$, $\rho\geq n-3(\delta-1)$. From ($10$), we have
\begin{center}
$\rho+1<n-\frac{n-2}{3}<n-\frac{\delta^3+\delta-1}{3}<n-3(\delta-1)\leq\rho$,
\end{center}
a contradiction.

Suppose $s\geq 3$, let $\varphi(\delta)=\delta^3+\delta-s-(\delta+1-s)(s+1)^2$ for $\delta>s$. Hence $\varphi^{\prime}(\delta)=3\delta^2+1-(s+1)^2>3s^2+1-s^2-2s-1=2s(s-1)>0$. Thus, $\varphi(\delta)>\varphi(s)=s^3-(s+1)^2>0$ for $\delta>s\geq3$. By $(11)$, we have
\begin{center}
$\rho+1<n-(\delta+1-s)(s+1)\leq\rho$,
\end{center}
a contradiction.

Thus $\rho_{\alpha}(G_s^2)< n-(\delta+1-s)(s+1)$, it follows that
\begin{center}
$\rho_{\alpha}(G_s^2)< n-(\delta+1-s)(s+1)=n-\delta-2-\Big((s-1)(\delta-s)+\delta-s-1\Big)\leq n-\delta-2$.
\end{center}

Since $K_{n-\delta-1}$ is a proper subgraph of $H$, we conclude
\begin{align}
\rho_{\alpha}(H)>\rho_{\alpha}(K_{n-\delta-1})=n-\delta-2
\end{align}
by Lemma \ref{le:2.1.}. It follows from ($3$), ($4$) and ($12$) that
\begin{align*}
\rho_{\alpha}(G)\leq\rho_{\alpha}(G_s^1)\leq\rho_{\alpha}(G_s^2)<n-\delta-2<\rho_{\alpha}(H),
\end{align*}
a contradiction.\\

\textbf{Case 3.} $s\geq\delta+1$.

Let $G_s^3=K_s\vee(K_{n-2s-1}\cup(s+1)K_1)$, we have $\rho_\alpha(G_s^1)\leq\rho_\alpha(G_s^3)$ by Lemma \ref{le:2.5.}. Therefore, it follows that
\begin{align}
\rho_{\alpha}(G)\leq\rho_{\alpha}(G_s^1)\leq\rho_{\alpha}(G_s^3).
\end{align}

We partition the vertex set of $G_s^3$ as $V(G_s^3)=V(K_s)\cup V((s+1)K_1)\cup V(K_{n-2s-1})$, where $V(K_s)=\{v_1,v_2,\ldots,v_s\}$, $V((s+1)K_1)=\{u_1,u_2,\ldots,u_{s+1}\}$ and $V(K_{n-2s-1})=\{w_1,w_2,\ldots,w_{n-2s-1}\}$. Clearly,
\begin{align*}
H=&G_s^3+\{u_iw_j|\delta+2\leq i\leq s+1, 1\leq j\leq n-2s-1\}
\\&+\{u_iu_j|\delta+2\leq i\leq s, i+1\leq j\leq s+1\}
\\&-\{v_iu_j|\delta+1\leq i\leq s, 1\leq j\leq \delta+1\}.
\end{align*}
Let $\rho^\prime=\rho_\alpha(G_s^3)$, $\rho^{\prime\prime}=\rho_\alpha(H)$, \textbf{y} be the Perron vector of $A_\alpha(G_s^3)$ with respect to $\rho^{\prime}$. By symmetry, \textbf{y} takes the same value (say $y_1$, $y_2$ and $y_3$) on the vertices of $V(K_s)$, $V((s+1)K_1)$ and $V(K_{n-2s-1})$, respectively. Then, by $A_\alpha(G_s^3)\textbf{y}=\rho^\prime\textbf{y}$, we have

$$\left\{
\begin{aligned}
&\rho^\prime y_{2}=(1-\alpha)s y_{1}+\alpha s y_{2},\\
&\rho^\prime y_{3}=(1-\alpha)s y_{1}+[n-(2-\alpha)s-2]y_{3} .
\end{aligned}
\right.$$

Since $\rho^\prime>0$ and $n\geq s+c(G-S)\geq s+s+2=2s+2$, we have
\begin{align*}
(\rho^\prime-\alpha s)(y_3-y_2)=(n-2s-2)y_3\geq 0,
\end{align*}
that is, $y_3\geq y_2$.

Similarly, let \textbf{z} be the Perron vector of $A_\alpha(H)$ with respect to $\rho^{\prime\prime}$. By symmetry, \textbf{z} takes the same value (say $z_1$, $z_2$ and $z_3$) on the vertices of $V(K_\delta)$, $V((\delta+1)K_1)$ and $V(K_{n-2\delta-1})$, respectively. Then, by $A_\alpha(H)\textbf{z}=\rho^{\prime\prime}\textbf{z}$, we have

$$\left\{
\begin{aligned}
&\rho^{\prime\prime} z_{2}=(1-\alpha)\delta z_{1}+\alpha \delta z_{2},\\
&\rho^{\prime\prime} z_{3}=(1-\alpha)\delta z_{1}+[n-(2-\alpha)\delta-2]z_{3},
\end{aligned}
\right.$$
which leads to
\begin{align*}
z_3=\frac{\rho^{\prime\prime}-\alpha\delta}{\rho^{\prime\prime}-\alpha\delta-(n-2\delta-2)}z_2.
\end{align*}

Since $G_s^3$ and $H$ are not regular, it follows that $\rho^\prime<n-1$ and $\rho^{\prime\prime}<n-1$. Recall that $n\geq 2s+2$, then $\delta+1\leq s\leq\frac{n-2}{2}$.

Suppose to the contrary that $\rho^\prime\geq\rho^{\prime\prime}$. Consider

\begin{align}
\textbf{z}^{T}(\rho^{\prime\prime}-\rho^\prime)\textbf{y}=&\textbf{z}^{T}(A_\alpha(H)-A_\alpha(G_s^3))\textbf{y}\nonumber
\\=&(1-\alpha)\textbf{z}^{T}(A(H)-A(G_s^3))\textbf{y}+\alpha\textbf{z}^{T}(D(H)-D(G_s^3))\textbf{y}\nonumber
\\=&(1-\alpha)\Big[\sum_{i=\delta+2}^{s+1}\sum_{j=1}^{n-2s-1}(y_{u_i}z_{w_j}+y_{w_j}z_{u_i})+ \sum_{i=\delta+2}^{s}\sum_{j=i+1}^{s+1}(y_{u_i}z_{u_j}+y_{u_j}z_{u_i})\nonumber
\\&-\sum_{i=\delta+1}^{s}\sum_{j=1}^{\delta+1}(y_{v_i}z_{u_j}+y_{u_j}z_{v_i})\Big]+\alpha\textbf{z}^{T}(D(H)-D(G_s^3))\textbf{y}\nonumber
\\=& (1-\alpha)(s-\delta)\Big[(n-2s-1)(y_{2}z_{3}+y_{3}z_{3})+(s-\delta-1)y_{2}z_{3}-(\delta+1)(y_{1}z_{2}+y_{2}z_{3})\Big]\nonumber
\\&+\alpha\Big[-(\delta+1)y_{1}z_{3}+(s-\delta)(y_{3}z_{3}-y_{2}z_{2})+(n-\delta-s-2)y_{2}z_{3}\Big]\nonumber
\\ > &(1-\alpha)(s-\delta)\Big[(n-2\delta-s-3)y_{2}z_{3}+(n-2s-1)y_{3}z_{3}-(\delta+1)y_{1}z_{2}\Big]-\alpha(\delta+1)y_{1}z_{3}\nonumber
\\ \geq &(1-\alpha)(s-\delta)\Big[(n-2\delta-s-2)y_{2}z_{3}-(\delta+1)y_{1}z_{2}\Big]-\alpha(\delta+1)y_{1}z_{3}\nonumber
\\ \geq &(1-\alpha)(s-\delta)\Big[(n-2\delta-\frac{n-2}{2}-2)y_{2}z_{3}-(\delta+1)y_{1}z_{2}\Big]-\alpha(\delta+1)y_{1}z_{3}\nonumber
\\=&(1-\alpha)(s-\delta)\Big[(\frac{n}{2}-2\delta-1)\frac{(1-\alpha)s}{\rho^\prime-\alpha s}y_{1}\frac{\rho^{\prime\prime}-\alpha\delta}{\rho^{\prime\prime}-\alpha\delta-(n-2\delta-2)}z_{2}-(\delta+1)y_{1}z_{2}\Big]\nonumber
\\&-\alpha(\delta+1)\frac{\rho^{\prime\prime}-\alpha\delta}{\rho^{\prime\prime}-\alpha\delta-(n-2\delta-2)}y_{1}z_{2}\nonumber
\\ \geq & \frac{(\delta+1)y_{1}z_{2}}{(\rho^\prime-\alpha s)(\rho^{\prime\prime}-\alpha\delta-(n-2\delta-2))}\Big[(1-\alpha)\Big((\rho^{\prime\prime}-\alpha\delta)(1-\alpha)(\frac{n}{2}-2\delta-1)\nonumber
\\ &-(\rho^\prime-\alpha s)(\rho^{\prime\prime}-\alpha\delta-(n-2\delta-2))\Big)-\alpha(\rho^{\prime\prime}-\alpha\delta)(\rho^\prime-\alpha s)\Big]\nonumber
\\>& \frac{(\delta+1)y_{1}z_{2}}{(\rho^\prime-\alpha s)(\rho^{\prime\prime}-\alpha\delta-(n-2\delta-2))}\cdot\nonumber
\\ &\Big((\rho^{\prime\prime}-\alpha\delta)(1-\alpha)^2(\frac{n}{2}-2\delta-1)+(\rho^{\prime\prime}-\alpha\delta)(\rho^\prime-\alpha s)\Big)\nonumber
\\>&\frac{(\delta+1)(\rho^{\prime\prime}-\alpha\delta)y_{1}z_{2}}{(\rho^\prime-\alpha s)(\rho^{\prime\prime}-\alpha\delta-(n-2\delta-2))}\Big((1-\alpha)^2(\frac{n}{2}-2\delta-1)+(n-2s-2)\Big)\nonumber
\\ \geq &0. \nonumber
\end{align}

It follows that $\rho^{\prime\prime}>\rho^\prime$ as $\textbf{z}^{T}\textbf{y}>0$, which contradicts the assumption that $\rho^\prime\geq\rho^{\prime\prime}$. Therefore $\rho^{\prime\prime}>\rho^\prime$. It follows from $(13)$ that
\begin{align*}
\rho_{\alpha}(G)\leq\rho_{\alpha}(G_s^1)\leq\rho_{\alpha}(G_s^3)=\rho^{\prime} < \rho^{\prime\prime}=\rho_{\alpha}(H).
\end{align*}

This completes the proof.
$\hfill\square$\\


\section{Proof of Theorem \ref{th:1.2.} }
\hspace{1.3em}

In this section, we give the proof of Theorem \ref{th:1.2.} and the following lemma are used further.

\noindent\begin{lemma}\label{le:3.1.}\cite{ABGD}
Let $G$ be a graph of order $n$ with $m$ edges, with no isolated vertices and let $\alpha\in[\frac{1}{2},1]$. Then

$$\rho_\alpha(G)\leq\frac{2m(1-\alpha)}{n-1}+\alpha n-1.$$
If $\alpha\in(\frac{1}{2},1)$ and $G$ is connected, equality holds if and only if $G=K_n$.
\end{lemma}

Now, we shall give a proof of Theorem\ref{th:1.2.}.\\

\noindent\textbf{Proof of Theorem\ref{th:1.2.}.}

Suppose to the contrary that $G$ is not a $\tau$-tough graph with $\tau\geq 2$ or $\frac{1}{\tau}$ is a positive integer. Then there exists some nonempty subset $S\subseteq V(G)$ such that $|S|<\tau(c(G-S)-1)$. Let $|S|=s$ and $c(G-S)=c$.\\

$(1)$ Since $\tau\geq 2$ be an integer, and $s\leq\tau(c-1)-1$, $G$ is a spanning subgraph of $G_1=K_{\tau(c-1)-1}\vee(K_{n_1}\cup K_{n_2}\cup\cdots\cup K_{n_c})$ for some integers $n_1\geq n_2\geq\cdots\geq n_c$ with $\sum\limits_{i=1}^{c}n_{i}=n-\tau(c-1)+1$. Combining this with Lemma \ref{le:2.2.}, we have

\begin{align}
\rho_\alpha(G)\leq\rho_\alpha(G_1)
\end{align}
with equality if and only if  $G\cong G_1$.

Let $G_2=K_{\tau(c-1)-1}\vee(K_{n-(\tau+1)(c-1)+1}\cup(c-1)K_{1})$. By Lemma \ref{le:2.5.}, we have
\begin{align}
\rho_\alpha(G_1)\leq\rho_\alpha(G_2)
\end{align}
with equality if and only if  $G_1\cong G_2$.\\

\textbf {Case 1}. $c=2$.

Then we have $G_2=K_{\tau-1}\vee(K_{n-\tau}\cup K_{1})$. By $(14)$ and $(15)$, we have
\begin{center}
$\rho_\alpha(G)\leq\rho_\alpha(K_{\tau-1}\vee(K_{n-\tau}\cup K_{1}))$
\end{center}
with equality if and only if  $G\cong K_{\tau-1}\vee(K_{n-\tau}\cup K_{1})$.\\

\textbf {Case 2}. $c\geq3$.

According to Lemma \ref{le:3.1.}, we have
\begin{align}
\rho_\alpha(G_2)\leq&\frac{2m(1-\alpha)}{n-1}+\alpha n-1\nonumber
\\=&\frac{(1-\alpha)\big[(n-c+1)(n-c)+2(\tau(c-1)-1)(c-1)\big]}{n-1}+\alpha n-1\nonumber
\\=&\frac{1-\alpha}{n-1}\big[(2\tau+1)c^2-(2n+4\tau+3)c+n^2+n+2\tau+2\big]+\alpha n-1.
\end{align}

Let $f(c)=(2\tau+1)c^2-(2n+4\tau+3)c+n^2+n+2\tau+2$. Observe that $n\geq(\tau+1)(c-1)$. Thus $3\leq c\leq\frac{n}{\tau+1}+1$. By a calculation,
\begin{center}
$f(3)-f(\frac{n}{\tau+1}+1)=\frac{(n-2\tau-2)(n-4\tau^2-5\tau-1)}{(\tau+1)^2}\geq 0$,
\end{center}
where the inequality follows from the fact that $n\geq4\tau^2+5\tau+1$. This implies that $f(c)$ attains its maximum value at $c=3$ when $3\leq c\leq\frac{n}{\tau+1}+1$. Combining this with $(16)$ and $n\geq\frac{8\tau(1-\alpha)-2\alpha+1}{3-4\alpha}$, we deduce that
\begin{align*}
\rho_\alpha(G_2)\leq &\frac{(1-\alpha)f(3)}{n-1}+\alpha n-1
\\=&\frac{(1-\alpha)(n^2-5n+8\tau+2))}{n-1}+\alpha n-1
\\=&n-2+\frac{-(3-4\alpha)n+8\tau(1-\alpha)-2\alpha+1}{n-1}
\\\leq& n-2+\frac{-(3-4\alpha)\frac{8\tau(1-\alpha)-2\alpha+1}{3-4\alpha}+8\tau(1-\alpha)-2\alpha+1}{n-1}
\\=&n-2.
\end{align*}

Since $K_{n-1}$ is a proper subgraph of $K_{\tau-1}\vee(K_{n-\tau}\cup K_1)$, we have
\begin{align}
\rho_{\alpha}(K_{\tau-1}\vee(K_{n-\tau}\cup K_1))>\rho_{\alpha}(K_{n-1})=n-2.
\end{align}

It follows from $(14), (15)$ and $(17)$ that
\begin{center}
$\rho_{\alpha}(G)\leq\rho_{\alpha}(G_1)\leq\rho_{\alpha}(G_2)\leq n-2=\rho_{\alpha}(K_{n-1})\leq\rho_{\alpha}(K_{\tau-1}\vee(K_{n-\tau}\cup K_1))$.
\end{center}

In conclusion, we have
\begin{center}
$\rho_{\alpha}(G)\leq\rho_{\alpha}(K_{\tau-1}\vee(K_{n-\tau}\cup K_1))$
\end{center}
with equality if and only if $G=K_{\tau-1}\vee(K_{n-\tau}\cup K_1)$, which is a contradiction to the $A_\alpha$-spectral radius condition of Theorem 1.2.\\

$(2)$ Notice that $\frac{1}{\tau}\geq 1$ is a positive number. We may assume that $b=\frac{1}{\tau}(b\geq1)$, and hence $c\geq bs+2$.
We know $G$ is a  spanning subgraph of $G^{\prime}=K_{s}\vee(K_{n_{1}}\cup K_{n_{2}}\cup\cdots\cup K_{n_{bs+2}})$ for some integers $n_{1}\geq n_{2}\geq\cdots\geq n_{bs+2}$ and $\sum\limits_{i=1}^{bs+2}n_{i}=n-s$.
By Lemma \ref{le:2.2.}, we obtain
\begin{align}
\rho_{\alpha}(G)\leq\rho_{\alpha}(G^{\prime})
\end{align}
with equality if and only if $G\cong G^{\prime}$. Let $G^{\prime\prime}=K_{s}\vee(K_{n-(b+1)s-1}\cup(b+1)K_{1})$.
 According to Lemma \ref{le:2.5.},  we have
\begin{align}
\rho_{\alpha}(G^{\prime})\leq\rho_{\alpha}(G^{\prime\prime})
\end{align}
with equality if and only if $G^{\prime}\cong G^{\prime\prime}$. It is obvious that $S\neq \varnothing$ since $S$ is a cut set, which means $s\geq 1$. Next we consider the following two cases based on the value of $s$.\\

\textbf{Case 1.} $s=1$.

In this case, $G^{\prime\prime}=K_{1}\vee(K_{n-b-2}\cup(b+1)K_{1})$. By $(18)$ and $(19)$, we obtain
 \begin{align}
 \rho_{\alpha}(G)\leq\rho_{\alpha}(G^{\prime\prime})
 \end{align}
with equality if and only if $G\cong G^{\prime\prime}$. By the assumption $\rho_{\alpha}(G^{\prime\prime})\leq\rho_{\alpha}(G)$, we have $\rho_{\alpha}(G^{\prime\prime})=\rho_{\alpha}(G)$ which implies $G\cong G^{\prime\prime}$.

Take $S=V(K_{1})$, thus $\tau(G^{\prime\prime})=\frac{s}{c(G^{\prime\prime}-S)-1}=\frac{1}{b+1}<\tau$, which implies $G^{\prime\prime}$ is not $\tau$-tough.
So $G\cong G^{\prime\prime}$.\\

\textbf{Case 2.} $s\geq2$.

According to Lemma \ref{le:3.1.}, we have
\begin{align}
\rho_\alpha(G^{\prime\prime})\leq&\frac{2m(1-\alpha)}{n-1}+\alpha n-1\nonumber
\\=&\frac{(1-\alpha)((b+2)bs^2-(2bn-3b-2)s+n^2-3n+2)}{n-1}+\alpha n-1\nonumber
\\=&\frac{1-\alpha}{n-1}[(b+2)bs^2-(2bn-3b-2)s+n^2-3n+2]+\alpha n-1.
\end{align}

Let $g(s)=(b+2)bs^2-(2bn-3b-2)s+n^2-3n+2$. It is clear that $n-(b+1)s-1\geq 1$. Thus $2\leq s\leq\frac{n-2}{b+1}$. By a calculation,

\begin{center}
$g(2)-g(\frac{n-2}{b+1})=\frac{(n-2b-4)((n-7)b^2-2b^3-5b-2)}{(b+1)^2}\geq 0$,
\end{center}
where the inequality follows from the fact that $n\geq2\tau^2+5\tau+\frac{2}{\tau}+8$. This implies that $g(s)$ attains its maximum value at $s=2$ when $2\leq s\leq\frac{n-2}{b+1}$. Combining this with $(21)$ and $n\geq\frac{(5-6\alpha)\tau^2+(13-14\alpha)\tau+4(1-\alpha)}{(1-2\alpha)\tau^2+(3-4\alpha)\tau}$, we deduce that

\begin{align}
\rho_\alpha(G^{\prime\prime})\leq &\frac{(1-\alpha)g(2)}{n-1}+\alpha n-1\nonumber
\\=&\frac{(1-\alpha)(n^2-(3+4b)n+4b^2+14b+6)}{n-1}+\alpha n-1\nonumber
\\=&n-b-2+\frac{-(3b-(2+4b)\alpha+1)n+(1-\alpha)(4b^2+14b+6)-b-1}{n-1}\nonumber
\\ \leq & n-b-2\nonumber
\\&+\frac{-(3b-(2+4b)\alpha+1)\frac{(5-6\alpha)\tau^2+(13-14\alpha)\tau+4(1-\alpha)}{(1-2\alpha)\tau^2+(3-4\alpha)\tau}+(1-\alpha)(4b^2+14b+6)-b-1}{n-1}\nonumber
\\=&n-b-2.
\end{align}

Since $K_{n-b-1}$ is a proper subgraph of $K_{1}\vee(K_{n-b-2}\cup (\frac{1}{\tau}+1)K_1)$, we have
\begin{align}
\rho_{\alpha}(K_{1}\vee(K_{n-b-2}\cup (\frac{1}{\tau}+1)K_1))>\rho_{\alpha}(K_{n-b-1})=n-b-2.
\end{align}

It follows from $(18), (19), (22)$ and $(23)$ that
\begin{center}
$\rho_{\alpha}(G)\leq\rho_{\alpha}(G^\prime)\leq\rho_{\alpha}(G^{\prime\prime})\leq n-b-2=\rho_{\alpha}(K_{n-b-1})<\rho_{\alpha}(K_{1}\vee(K_{n-b-2}\cup(b+1) K_1))$.
\end{center}

In conclusion, we have
\begin{center}
$\rho_{\alpha}(G)<\rho_{\alpha}(K_{1}\vee(K_{n-\frac{1}{\tau}-2}\cup (\frac{1}{\tau}+1)K_1))$.
\end{center}

This completes the proof.
$\hfill\square$\\

\textbf{Declaration of competing interest}\\

The authors declare that they have no known competing financial interests or personal relationships that could have appeared to influence the work reported in this paper.\\

\textbf{Date availability}\\

No date was used for the research described in the article.

\end{document}